\documentclass[11pt]{amsart}
\usepackage{amssymb,amscd,amsmath,amsthm}

\theoremstyle{plain}
\newtheorem{thm}{Theorem}[section]

\newtheorem{lem}[thm]{Lemma}

\theoremstyle{definition}

\newtheorem{ex}[thm]{Example}

\theoremstyle{remark}
\newtheorem{rem}{Remark}

\newtheorem*{remark*}{Remark}

\numberwithin{equation}{section}

\begin{document}
\title{Integer points in domains and adiabatic limits}

\author[Y. A. Kordyukov]{Yuri A. Kordyukov}
\address{Institute of Mathematics\\
         Russian Academy of Sciences\\
         112~Chernyshevsky str.\\ 450008 Ufa\\ Russia} \email{yurikor@matem.anrb.ru}

\author[A. A. Yakovlev]{Andrey A. Yakovlev}
\address{Department of Mathematics\\
Ufa State Aviation Technical University\\ 12 K. Marx str.\\ 450000
Ufa\\ Russia} \email{yakovlevandrey@yandex.ru}
\thanks{Supported by the Russian Foundation of Basic Research
(grant no. 09-01-00389)}



\begin{abstract}
We prove an asymptotic formula for the number of integer points in a
family of bounded domains with smooth boundary in the Euclidean
space, which remain unchanged along some linear subspace and expand
in the directions, orthogonal to this subspace. A more precise
estimate for the remainder is obtained in the case when the domains
are strictly convex.

Using these results, we improved the remainder estimate in the
adiabatic limit formula (due to the first author) for the eigenvalue
distribution function of the Laplace operator associated with a
bundle-like metric on a compact manifold equipped with a Riemannian
foliation in the particular case when the foliation is a linear
foliation on the torus and the metric is the standard Euclidean
metric on the torus.
\end{abstract}

\date{\today}

\maketitle

\section{Statement of the problem and the main results}
A classical problem on integer points distribution consists in the
study of the asymptotic behavior of the number of points of the
integer lattice ${\mathbb Z}^n$ in a family of homothetic domains in
${\mathbb R}^n$. This problem is originated in the Gauss problem on
the number of integer points in the disk, where it is directly
related with the arithmetic problem on the number of representations
of an integer as a sum of two squares, and sufficiently well studied
(see, for instance, books
\cite{Fricker,Gruber-Lekker,Huxley,Kraetzel} and the references
therein).

In this paper we investigate much less studied problem on counting
of integer points in a family of anisotropically expanding domains.
More precisely, let $F$ be a $p$-dimensional linear subspace of
$\mathbb {R}^n$ and $H=F^{\bot}$ the $q$-dimensional orthogonal
complement of $F$ with respect to the standard inner product
$(\cdot,\cdot)$ in $\mathbb {R}^n$, $p+q=n$. For any
$\varepsilon>0$, consider the linear transformation $T_\varepsilon :
\mathbb {R}^n\to \mathbb {R}^n$ given as follows:
\[
T_\varepsilon(x)=\begin{cases} x, & \text{if}\ x\in F, \\
\varepsilon^{-1}x, & \text{if}\ x\in H.
\end{cases}
\]
For any bounded domain $S$ with smooth boundary in $\mathbb{R}^n$,
we put
\begin{equation}\label{e:nS}
n_\varepsilon (S)=\# (T_\varepsilon(S)\cap \mathbb{Z}^n), \quad
\varepsilon>0.
\end{equation}

The main goal of the paper is to study the asymptotic behavior of
the function $n_\varepsilon (S)$ as $\varepsilon\to 0$. Before we
state the main results of the paper, we introduce some auxiliary
notions.

Let $\Gamma = \mathbb{Z}^n\cap F$. $\Gamma$ is a free abelian group.
Denote by $r=\operatorname{rank} \Gamma \leq p$ the rank of
$\Gamma$. For $r \geq 1$, denote by $(\ell_1, \ell_2, ... \ell_r)$
some base in $\Gamma$. Let $V$ be the $r$-dimensional subspace of
$\mathbb R^n$ spanned by $(\ell_1, \ell_2, ... \ell_r)$. Observe
that $\Gamma$ is a lattice in $V$. Denote by $Q\subset V$ the
parallelepiped spanned by the vectors $(\ell_1, \ell_2, ... \ell_r)$
and by $|Q|$ its $r$-dimensional Euclidean volume:
\[
|Q|={\rm vol}_r(\ell_1, \ell_2, ... \ell_r)={\rm vol}(V/ \Gamma).
\]

Let $\Gamma^*$ denote the lattice in $V$, dual to the lattice
$\Gamma$:
\[
\Gamma^*=\{ \gamma^*\in V :(\gamma^*,\Gamma)\subset\mathbb Z\}.
\]
For $r=0$, the groups $\Gamma$ and $\Gamma^*$ are trivial, and it is
natural to put $|Q|=1$.

For any $x\in V$, we denote by $P_{x}$ the $(n-r)$-dimensional
affine subspace of $\mathbb {R}^n$, passing through $x$ orthogonal
to $V$.

\begin{thm}\label{mainthm1}
For any bounded open set $S$ with smooth boundary in $\mathbb{R}^n$,
the formula holds:
\begin{equation}\label{f:lattice_points1}
 n_\varepsilon (S)\\ =
\frac{\varepsilon^{-q}}{|Q|}
 \sum_{\gamma^*\in\Gamma^*} {\rm vol}_{n-r} (P_{\gamma^*}\cap S)
+O(\varepsilon^{\frac{1}{p-r+1}-q}), \quad \varepsilon \to 0.
\end{equation}
\end{thm}

Note that, in the general case, the leading term in the asymptotic
formula for $n_\varepsilon (S)$ as $\varepsilon\to 0$ was unknown.

\begin{rem}
In the case when $F$ is trivial, we have $p=r=0$, $q=n$. The problem
is reduced to the classical problem on the asymptotics of the number
of integer points in a family of homothetic domains in
$\mathbb{R}^n$. The formula \eqref{f:lattice_points1} is reduced to
the classical formula, going back to Gauss:
\[
 n_\varepsilon (S) =
\varepsilon^{-n}{\rm vol}_{n}(S) +O(\varepsilon^{1-n}), \quad
\varepsilon \to 0.
\]
\end{rem}

\begin{rem}\label{r:r=0}
In the case when $\Gamma$ is trivial, we have $r=0$. The formula
\eqref{f:lattice_points1} takes the form
\[
 n_\varepsilon (S)=\varepsilon^{-q}{\rm vol}_{n} (S)
+O(\varepsilon^{\frac{1}{p+1}-q}), \quad \varepsilon \to 0.
\]
\end{rem}

In order to obtain a more precise estimate for the remainder, we
need to impose some restrictions to the boundary of the domain $S$.
We give just one result of similar type.

\begin{thm}\label{mainthm}
For any bounded open set $S$ with smooth boundary in $\mathbb{R}^n$
such that, for any $x\in F$, the intersection $S\cap \{x+H\}$ is
strictly convex, the formula holds:
\begin{equation}\label{f:lattice_points}
 n_\varepsilon (S)\\ =
\frac{\varepsilon^{-q}}{|Q|}
 \sum_{\gamma^*\in\Gamma^*} {\rm vol}_{n-r} (P_{\gamma^*}\cap S)
+O(\varepsilon^{k-q}), \quad \varepsilon \to 0,
\end{equation}
where
\[
k=\begin{cases} \frac{q+1}{2(p-r+1)} & \text{if}\ \frac{q-1}{2}\leq
p-r\\ \frac{2q}{q+1+2(p-r)} & \text{if}\ \frac{q-1}{2}>p-r.
\end{cases}
\]
\end{thm}

\begin{rem}
In the case when $F$ is trivial, we get $k=2n/(n+1)$, and the
formula \eqref{f:lattice_points1} is reduced to the following
formula:
\[
 n_\varepsilon (S) =
\varepsilon^{-n}{\rm vol}_{n}(S)
+O(\varepsilon^{-n+2-\frac{2}{n+1}}), \quad \varepsilon \to 0.
\]
This formula was proved by Randol \cite{Randol66I,Randol66II}.
\end{rem}

\begin{ex} As the simplest non-trivial example, we consider the case when
$n=2$ and $p=1$. Thus, let $F$ be the one-dimensional linear
subspace of $\mathbb{R}^2$ spanned by $(1,\alpha)\in \mathbb{R}^2$.
Then its orthogonal complement $H$ is spanned by $(-\alpha,1)\in
\mathbb{R}^2$. One can distinguish two cases.
\medskip
\par
{\it Case 1: $\alpha\not\in\mathbb{Q}$.} In this case, $\Gamma$ is
trivial, therefore, $r=0$. Moreover, $\Gamma^*$ is trivial, and, as
mentioned above, it is natural to put $|Q|=1$. For any bounded
domain $S$ with smooth boundary in ${\mathbb R}^2$, we get
\[
n_{\varepsilon}(S)=\varepsilon^{-1} {\rm
area}(S)+O(\varepsilon^{-1/2}), \quad \varepsilon\rightarrow 0.
\]

{\it Case 2: $\alpha\in\mathbb{Q}$.} Write $\alpha=\frac{p}{q}$,
where $p$ and $q$ are coprime numbers. Then $\Gamma$ is generated by
$\ell_1=(q,p)$. Therefore,
\[
|Q|=|\ell_1|=\sqrt{p^2+q^2}.
\]
The dual lattice $\Gamma^*$ is generated by
$\frac{1}{p^2+q^2}\ell_1$. An arbitrary element of $\Gamma^*$ has
the form
\[
\gamma^*=\frac{1}{p^2+q^2}k \ell_1, \quad k\in {\mathbb Z}.
\]
The corresponding subspace $P_{\gamma^*}$ is the line $L_k$ on the
plane ${\mathbb R}^2$ given by the equation $qx+py-k=0$. Therefore,
for any bounded domain $S$ with smooth boundary in ${\mathbb R}^2$,
we get
\[
n_{\varepsilon}(S) = \varepsilon^{-1} \frac{1}{\sqrt{p^2+q^2}}
\sum_{k\in {\mathbb Z}} |S\cap L_k|+O(1), \quad
\varepsilon\rightarrow 0.
\]
\end{ex}

The problem on counting of integer points in anisotropically
expanding domains in a slightly different context was studied in
considerable detail in
\cite{Skriganov89,Skriganov94,Nikichine-Skriganov95,Nikichine-Skriganov98}.
More precisely, these papers were devoted to the estimates of the
number $N(S,\Gamma)$ of points of a lattice $\Gamma$ in ${\mathbb
R}^n$ lying in a bounded domain $S$ for special domains and
lattices. Let us briefly describe some results of these papers and
show how they can be applied to the problem under consideration.

Let $k$ be a totally real algebraic number field of degree $n$,
$\sigma$ the canonical embedding of $k$ in the Euclidean space
${\mathbb R}^n$, $M\subset k$ an arbitrary ${\mathbb Z}$-module of
rank $n$ and $\Gamma_M=\sigma(M)$ the corresponding algebraic
lattice in ${\mathbb R}^n$. Let $\Pi\subset {\mathbb R}^n$ be an
$n$-dimensional parallelepiped centered at the origin and with edges
parallel to the coordinated axes:
\[
\Pi=\prod_{j=1}^n(-a_j,a_j).
\]

For $\lambda=(\lambda_1,\ldots,\lambda_n)\in {\mathbb R}^n$, we put
\[
\operatorname{Nm} \lambda=\prod_{j=1}^n\lambda_j.
\]
Let $\{e_1,\ldots,e_n\}$ be the standard basis of ${\mathbb R}^n$.
Any vector $\lambda=(\lambda_1,\ldots,\lambda_n)\in {\mathbb R}^n$
defines a linear transformation of ${\mathbb R}^n$ by
\[
\lambda\cdot e_j=\lambda_je_j,\quad j\in \{1,\ldots,n\}.
\]

By \cite[Theorem 1.1]{Skriganov94} (see also
\cite{Skriganov89,Nikichine-Skriganov98}), for any $\lambda\in
{\mathbb R}^n$ we have the estimate
\begin{equation}\label{e:Skr}
\left|N(\lambda\cdot \Pi,\Gamma_M)-|\operatorname{Nm}
\lambda|\frac{{\rm vol}_{n}(\Pi)}{d(\Gamma_M)}\right|< C[\ln
(2+|\operatorname{Nm} \lambda|)]^{n-1},
\end{equation}
where $d(\Gamma_M)$ is the volume of a fundamental domain of the
lattice $\Gamma_M$ and $C>0$ is a constant, independent of $\lambda
$.

In \cite{Nikichine-Skriganov95,Nikichine-Skriganov98}, this result
was extended to the case of algebraic lattices associated with an
arbitrary algebraic number field as well as to a wider class of
domains. As an application, the authors obtain nontrivial remainder
estimates in the asymptotic formula due to Lang \cite{Lang} for the
number of elements of an algebraic number field in a parallelotope
determined by the canonical system of valuations.

The study carried out in
\cite{Skriganov89,Skriganov94,Nikichine-Skriganov95,Nikichine-Skriganov98}
shows that the remainder estimate in the asymptotic formula for
$N(\lambda \cdot S,\Gamma)$ turns out to be very sensitive to the
Diophantian properties of the lattice $\Gamma$ and to the geometric
properties of the domain $S$.

In order to apply the results described above to the problems
considered in this paper, let us fix $p\in \{1,\ldots,n\}$ and make
use of the formula \eqref{e:Skr} for
$\lambda=\lambda(\varepsilon)\in {\mathbb R}^n$ of the form
\[
\lambda_j=\begin{cases} 1, & \text{if}\ 1\leq j\leq p,\\
\varepsilon^{-1}, & \text{if}\ p+1\leq j\leq n,\end{cases}
\]
for some $\varepsilon>0$. We get
\[
\left|N(\lambda(\varepsilon) \cdot
\Pi,\Gamma_M)-\varepsilon^{-q}\frac{{\rm
vol}_{n}(\Pi)}{d(\Gamma_M)}\right|< C_1|\ln \varepsilon|^{n-1},
\quad \varepsilon \to 0,
\]
where $C_1>0$ is a constant, independent of $\varepsilon $, and
$p+q=n$.

On the other hand, the lattice $\Gamma_M$ can be represented as
$\Gamma_M=A({\mathbb Z}^n)$ with some linear isomorphism $A$ of
${\mathbb R}^n$. Therefore, we have
\[
N(\Pi,\Gamma_M)=N(\Pi, A({\mathbb Z}^n))=N(A^{-1}(\Pi),{\mathbb
Z}^n).
\]
It is not difficult to see from here that
\[
N(\lambda(\varepsilon)\cdot \Pi,\Gamma_M)
=n_\varepsilon(A^{-1}(\Pi)),
\]
where the right-hand side of this identity is defined by
\eqref{e:nS} with the subspace $F\subset {\mathbb R}^n$ spanned by
$A^{-1}(e_1),\ldots,A^{-1}(e_p)$ and the subspace $H\subset {\mathbb
R}^n$ spanned by $A^{-1}(e_{p+1}),\ldots,A^{-1}(e_n)$. Observe that
${\mathbb Z}^n\cap F=\{0\}$.

Thus, finally we get
\[
n_\varepsilon(A^{-1}(\Pi))=\varepsilon^{-q}{\rm
vol}_{n}(A^{-1}(\Pi))+O(|\ln \varepsilon|^{n-1}), \quad \varepsilon
\to 0,
\]
which is the formula \eqref{f:lattice_points1} with a more precise
remainder estimate in the case when $S= A^{-1}\Pi$ is a
parallelepiped centered at the origin and with edges parallel to the
vectors $A^{-1}(e_{1}),\ldots,A^{-1}(e_n)$. Note that, in this case,
the subspaces $F$ and $H$ have rather special form and, in general,
are not orthogonal.

Similarly, one can use the results of
\cite{Nikichine-Skriganov95,Nikichine-Skriganov98} and get more
precise remainder estimates in \eqref{f:lattice_points1} for some
special subspaces $F$ and $H$ and domains $S$.

It would be interesting to continue the study of remainder estimates
in the asymptotic formula~\eqref{f:lattice_points1}, depending on
geometry of a domain $S$ and properties of $F$ and $H$.

The authors are grateful to the referee for useful remarks.

\section{Applications to adiabatic limits}
It is well known that the Gauss problem on counting of integer
points in the disk is equivalent to the problem on the asymptotic
behavior of the eigenvalue distribution function of some elliptic
differential operator on a compact manifold, namely, of the Laplace
operator on a torus. In the case under consideration, there is also
an equivalent asymptotic spectral problem, namely, the problem on
the asymptotic behavior of the eigenvalue distribution function of
the Laplace operator on a torus in the adiabatic limit associated
with a linear foliation.

Consider the $n$-dimensional torus $\mathbb T^n=\mathbb
{R}^n/\mathbb{Z}^n$. Let $\mathcal F$ be a linear foliation on
$\mathbb T^n$: the leaf $L_x$ of $\mathcal F$ through $x\in \mathbb
T^n$ has the form:
\[
L_x=x+F \mod \mathbb{Z}^n.
\]
The decomposition of $\mathbb{R}^n$ into the direct sum of subspaces
$\mathbb{R}^n=F\bigoplus H$ induces the corresponding decomposition
$g=g_{F}+g_{H}$ of the Euclidean metric into the sum of the
tangential and transversal components. Define a one-parameter family
$g_{\varepsilon}$ of Euclidean metrics on $\mathbb{R}^n$ by
\begin{equation*}
g_{\varepsilon}=g_{F} + {\varepsilon}^{-2}g_{H}, \quad \varepsilon >
0.
\end{equation*}
We will also consider the metrics $g_{\varepsilon}$ as Riemannian
metrics on $\mathbb{T}^n$.

Let $A=(a_1,\ldots,a_n)\in \mathbb{R}^n$ be some point. Define a
$1$-form $\mathbf A$ on $\mathbb{T}^n$ by
\[
{\mathbf A}=\sum_{j=1}^n a_jdx_j.
\]
Consider the operator $d-2\pi i{\mathbf A}$, acting from
$C^\infty(\mathbb T^n)$ to the space $\Omega^1(\mathbb{T}^n)$ of
smooth $1$-forms on $\mathbb{T}^n$, where $d$ is the de Rham
differential, and $\mathbf A$ is the multiplication operator by
$\mathbf A$. Let $(d-2\pi i{\mathbf A})^{*}_{g_\varepsilon}:
\Omega^1(\mathbb{T}^n) \rightarrow C^\infty(\mathbb T^n)$ be the
adjoint of $d-2\pi i{\mathbf A}$ with respect to the inner products
in $C^\infty(\mathbb T^n)$ and $\Omega^1(\mathbb{T}^n)$ determined
by $g_{\varepsilon}$.

For any $\varepsilon>0$, consider the operator $H_\varepsilon$ in
$C^\infty(\mathbb T^n)$ defined by
\[
H_\varepsilon=(d-2\pi i{\mathbf A})^*_{g_\varepsilon}(d-2\pi
i{\mathbf A}).
\]
In the local coordinates $(x_1,x_2,\ldots,x_n)$ of the space
$\mathbb R^n$, the operator $H_\varepsilon$ is written in the form
\[
H_\varepsilon =\sum_{j,\ell=1}^n g_\varepsilon^{j\ell}
\left(\frac{\partial}{\partial x_j} -2\pi
ia_j\right)\left(\frac{\partial}{\partial x_\ell} - 2\pi
ia_\ell\right),
\]
where $g_\varepsilon^{j\ell}$ are the elements of the inverse matrix
of $g_\varepsilon $.

The operator $H_\varepsilon$ can be considered as the magnetic
Schr\"odinger operator on the torus $\mathbb T^n$, associated with
the metric $g_\varepsilon$, with the constant magnetic potential
${\mathbf A}$. It has a complete orthogonal systems of
eigenfunctions
\[
U_{k}(x)=e^{2\pi i (k,x)}, \quad x\in \mathbb R^n,\quad k\in \mathbb
Z^n,
\]
with the corresponding eigenvalues
\[
\lambda_{k}=(2\pi)^2
\|k-A\|^2_{g_\varepsilon^{-1}}=(2\pi)^2\sum_{j,\ell=1}^n
g_\varepsilon^{j\ell} \left(k_j -a_j\right) \left(k_\ell -
a_\ell\right).
\]

Denote by $N_\varepsilon(\lambda)$ the eigenvalue distribution
function of $H_\varepsilon$:
\[
N_\varepsilon(\lambda)=\sharp\{k\in \mathbb Z^n:\lambda_{k}<\lambda
\}, \quad \lambda \in {\mathbb R}.
\]
It is easy to see that
\[
n_\varepsilon (B_{\sqrt{\lambda}}(A))=N_\varepsilon (4\pi^2\lambda),
\quad \lambda \in {\mathbb R}.
\]
Thus, the problem on the asymptotic behavior of the number
$n_\varepsilon (B_{\sqrt{\lambda}}(A))$ of integer points in the
ellipsoid $T_\varepsilon(B_{\sqrt{\lambda}}(A))$ as $\varepsilon\to
0$ is equivalent to the problem on the asymptotic behavior of the
eigenvalue distribution function $N_\varepsilon(\lambda)$ as
$\varepsilon\to 0$. The limiting procedure $\varepsilon\to 0$ will
be called adiabatic limit. This notion was introduced by Witten in
1985 in the study of global anomalies in string theory. We refer the
reader to a survey paper \cite{bedlewo2-andrey} for some historic
remarks and references.

In~\cite{adiab} (see also \cite{asymp}), the first author computed
the leading term of the asymptotics of the eigenvalue distribution
function of the Laplace operator associated with a bundle-like
metric on a compact manifold equipped with a Riemannian foliation,
in adiabatic limit. The linear foliation on the torus is a
Riemannian foliation, and a Euclidean metric on the torus is
bundle-like. As a straightforward consequence of
Theorem~\ref{mainthm}, we obtain a more precise estimate for the
remainder in the asymptotic formula of~\cite{adiab} for this
particular case (see also \cite{torus-eng}).

\begin{thm}\label{mainthm2}
For $\lambda>0$, the following asymptotic formula holds as
$\varepsilon \to 0$:
\[
N_\varepsilon (\lambda) = \varepsilon^{-q} \frac{\omega_{n-r}}{|Q|}
 \sum_{\gamma^*\in\Gamma^*} \left(\frac{\lambda}{4\pi^2}-|\gamma^*-A|^2\right)^{(n-r)/2}
+O(\varepsilon^{k-q}),
\]
where $\omega_{n-r}$ is the volume of the unit ball in ${\mathbb
R}^{n-r}$ and
\[
k=\begin{cases} \frac{q+1}{2(p-r+1)}, & \text{if}\ \frac{q-1}{2}\leq
p-r\\ \frac{2q}{q+1+2(p-r)}, & \text{if}\ \frac{q-1}{2}>p-r.
\end{cases}
\]
\end{thm}

\section{Proofs of the main results}
We will use a method based on the Poisson summation formula and the
method of stationary phase (Van der Corput, Randol
\cite{Randol66I,Randol66II}, Colin de Verdi{\`e}re \cite{ColindeV}).
First of all, we observe that we have the inclusion
\begin{equation}\label{e:TS}
{\mathbb Z}^n \subset \bigsqcup_{\gamma^*\in \Gamma^*} P_{\gamma^*}.
\end{equation}
Indeed, let $k\in {\mathbb Z}^n$. Denote by $\pi_V : {\mathbb
R}^n\to V$ the orthogonal projection on $V$. Then, for any
$\gamma\in \Gamma$, we have
\[
(\pi_V(k),\gamma)=(k,\gamma)\in {\mathbb Z},
\]
since $\gamma\in {\mathbb Z}^n$ and $\gamma\in F\subset V$. Hence,
$\pi_V(k)\in \Gamma^*$, that immediately implies \eqref{e:TS}.

For any $\gamma^*\in \Gamma^*$, denote
\[
{\mathbb Z}^n_{\gamma^*}={\mathbb Z}^n\bigcap P_{\gamma^*}=\{ k\in
{\mathbb Z}^n : \pi_V(k)=\gamma^*\}.
\]
We identify the affine subspace $P_{\gamma^*}$ with the linear space
$V^\bot$, fixing an arbitrary point $k_{\gamma^*}\in {\mathbb
Z}^n_{\gamma^*}$:
\[
P_{\gamma^*}=k_{\gamma^*}+V^\bot.
\]
It is easy to see that
\[
{\mathbb Z}^n_{\gamma^*}=k_{\gamma^*}+\Gamma^\bot,
\]
where
\[
\Gamma^\bot={\mathbb Z}^n\bigcap V^\bot
\]
is a lattice in $V^\bot$.

Let us observe the following relation, which will be needed later:
\begin{equation}\label{e:vol}
{\rm vol}(V^\bot/\Gamma^\bot)=|Q|.
\end{equation}
For its proof, let us choose some base $(\ell_1, \ldots, \ell_r)$ in
$\Gamma$. Denote by $(\ell^*_1, \ldots, \ell^*_r)$ the dual base in
$\Gamma^*$: $(\ell_i,\ell^*_j)=\delta_{ij}$ for any
$i,j=1,\ldots,r$. For any $i=1,\ldots,r$, choose some $k^*_i\in
{\mathbb Z}^n$ such that $\pi_V(k^*_i)=\ell^*_i$. Let $(k^\bot_1,
\ldots, k^\bot_{n-r})$ be a base in $\Gamma^\bot$. Using
\eqref{e:TS}, it easy to show that the system $(k^*_1,\ldots, k^*_r,
k^\bot_1, \ldots, k^\bot_{n-r})$ is a base in ${\mathbb Z}^n$.
Therefore, for the volume of the parallelepiped spanned by
$(k^*_1,\ldots, k^*_r, k^\bot_1, \ldots, k^\bot_{n-r})$, one has the
following formula:
\[
{\rm vol}_n(k^*_1,\ldots, k^*_r, k^\bot_1, \ldots,
k^\bot_{n-r})={\rm vol}({\mathbb R}^n/ {\mathbb Z}^n)=1.
\]
On the other hand, since $k^*_i-\ell^*_i\in V^\bot$ for any
$i=1,\ldots,r$, and the systems $(\ell^*_1,\ldots, \ell^*_r)$ and
$(k^\bot_1, \ldots, k^\bot_{n-r})$ are mutually orthogonal, we get
\begin{align*}
{\rm vol}_n(k^*_1,\ldots, k^*_r, k^\bot_1, \ldots, k^\bot_{n-r})= &
{\rm vol}_n(\ell^*_1,\ldots, \ell^*_r, k^\bot_1, \ldots,
k^\bot_{n-r})\\ = & {\rm vol}_r(\ell^*_1,\ldots, \ell^*_r){\rm
vol}_{n-r}(k^\bot_1, \ldots, k^\bot_{n-r})\\ = & {\rm
vol}(V/\Gamma^*){\rm vol}(V^\bot/\Gamma^\bot).
\end{align*}
Thus, we have
\[
{\rm vol}(V^\bot/\Gamma^\bot)=\frac{1}{{\rm vol}(V/\Gamma^*)}.
\]
In order to complete the proof of \eqref{e:vol}, it remains to apply
the well known relation
\[
{\rm vol}(V/\Gamma){\rm vol}(V/\Gamma^*)=1.
\]

Thus, we can write
\begin{equation}\label{e:ne}
n_\varepsilon (S)=\sum_{\gamma^*\in \Gamma^*}n_\varepsilon
(S,{\gamma^*}),
\end{equation}
where
\[
n_\varepsilon (S,{\gamma^*})=\# (T_\varepsilon(S)\cap {\mathbb
Z}^n_{\gamma^*}).
\]
Note that, since $S$ is bounded, the sum in the right hand side
of~\eqref{e:ne} has finitely many non-vanishing terms.

Fix $\gamma^*\in \Gamma^*$. Let $\chi_{S_{\gamma^*}}$ be the
indicator of the set $S_{\gamma^*}=S\bigcap P_{\gamma^*}$. It is
easy to see that
\begin{align*}
n_\varepsilon (S,{\gamma^*})& =\sum_{k\in {\mathbb
Z}^n_{\gamma^*}}\chi_{T_\varepsilon(S_{\gamma^*})}(k)=\sum_{\gamma\in
\Gamma^\bot}\chi_{T_\varepsilon(S_{\gamma^*})}(k_{\gamma^*}+\gamma)\\
& = \sum_{\gamma\in \Gamma^\bot} \chi_{S_{\gamma^*}}
(T_{\varepsilon^{-1}}(k_{\gamma^*}+\gamma))\\
& = \sum_{\gamma\in \Gamma^\bot} \chi_{S_{\gamma^*}}
(k_{\gamma^*}+(T_{\varepsilon^{-1}}(k_{\gamma^*})-k_{\gamma^*})+T_{\varepsilon^{-1}}(\gamma))).
\end{align*}

The space $V^\bot$ decomposes into the direct sum
\begin{equation}\label{e:Vbotdecomp}
V^\bot=F_V\bigoplus H,
\end{equation}
where $F_V=F\cap V^\bot$. We will write the decomposition of $x\in
V^\bot$, corresponding to \eqref{e:Vbotdecomp}, as follows:
\[
x=x_F+x_H, \quad x_F\in F_V, x_H\in H.
\]
Note that
\[
T_{\varepsilon}(x)=x_F +\varepsilon^{-1}x_H.
\]

Let $\rho \in C^\infty_0({\mathbb R})$ be an even function such that
$0\leq \rho(x) \leq 1$ for any $x\in {\mathbb R}$ and ${\rm
supp}\,\rho\subset (-1,1)$. For any $t_F>0$ and $t_H>0$, define a
function $\rho_{t_F,t_H} \in C^\infty_0(V^\bot)$ by
\[
\rho_{t_F,t_H}(x)=\frac{c}{t_F^{p-r}t_H^{q}}\rho
\left(\left(t^{-2}_Fx^2_F+t^{-2}_Hx^2_H\right)^{1/2}\right), \quad
x\in V^\bot,
\]
where the constant $c>0$ is chosen so that $\int_{V^\bot}
\rho_{1,1}(x)\,dx=1$. The function $\rho_{t_F,t_H}$ is supported in
the ellipsoid
\[
B(0,t_F,t_H)=\left\{x\in V^\bot :
\frac{x^2_F}{t^2_F}+\frac{x^2_H}{t^2_H}<1\right\}.
\]

Define the function $n_{\varepsilon,t_F,t_H}(S,\gamma^*)$ by
\[
n_{\varepsilon,t_F,t_H}(S,\gamma^*)=\sum_{k\in {\mathbb
Z}^n_{\gamma^*}}(\chi_{T_\varepsilon(S_{\gamma^*})}\ast
\rho_{t_F,t_H})(k),
\]
where the function $\chi_{T_\varepsilon(S_{\gamma^*})} \ast
\rho_{t_F,t_H}\in C^\infty_0(P_{\gamma^*})$ is defined by
\[
(\chi_{T_\varepsilon(S_{\gamma^*})} \ast \rho_{t_F,t_H})(y)
=\int_{V^\bot} \chi_{T_\varepsilon(S_{\gamma^*})} (y-x)
\rho_{t_F,t_H}(x)\, dx, \quad y\in P_{\gamma^*}.
\]

For any domain $S\subset P_{\gamma^*}$ and for any $t_F>0$ and
$t_H>0$, denote
\[
S_{t_F,t_H}=\bigcup_{x\in S}(x+B(0,t_F,t_H)),
\]
and
\[
S_{-t_F,-t_H} = P_{\gamma^*} \setminus (P_{\gamma^*} \setminus
S)_{t_F,t_H}.
\]
It is easy to see that, for any $\varepsilon>0$, $t_F>0$ and
$t_H>0$, one has
\[
T_\varepsilon(S_{t_F,\varepsilon t_H})=(T_\varepsilon(S))_{t_F,t_H}.
\]

\begin{lem}
For any $\varepsilon>0$, $t_F>0$ and $t_H>0$, the following
inequalities hold:
\[
n_{\varepsilon,t_F,t_H}((S_{\gamma^*})_{-t_F,-\varepsilon
t_H},\gamma^*) \leq n_\varepsilon (S,{\gamma^*}) \leq
n_{\varepsilon,t_F,t_H}((S_{\gamma^*})_{t_F,\varepsilon
t_H},\gamma^*).
\]
\end{lem}

\begin{proof}
Suppose that $k\in \mathbb Z^n\cap T_\varepsilon(S_{\gamma^*})$. For
any $x\in B(0,t_F,t_H)$, the point $k-x$ belongs to
$(T_\varepsilon(S_{\gamma^*}))_{t_F,t_H}$. Therefore,
$\chi_{T_\varepsilon((S_{\gamma^*})_{t_F,\varepsilon t_H})} (k-x)=1$
and
\[
(\chi_{T_\varepsilon((S_{\gamma^*})_{t_F,\varepsilon t_H})} \ast
\rho_{t_F,t_H})(k)=\int_{V^\bot} \rho_{t_F,t_H}(x) dx=1
=\chi_{T_\varepsilon(S_{\gamma^*})} (k).
\]
If $k\not\in \mathbb Z^n\cap T_\varepsilon(S_{\gamma^*})$, then
\[
\chi_{T_\varepsilon(S_{\gamma^*})}(k)=0\leq
(\chi_{T_\varepsilon((S_{\gamma^*})_{t_F,\varepsilon t_H})} \ast
\rho_{t_F,t_H})(k).
\]
We get
\begin{multline*}
n_\varepsilon (S,{\gamma^*})=\sum_{k\in \mathbb
Z^n_{\gamma^*}}\chi_{T_\varepsilon(S_{\gamma^*})}(k)\\ \leq
\sum_{k\in \mathbb Z^n_{\gamma^*}}
(\chi_{T_\varepsilon((S_{\gamma^*})_{t_F,\varepsilon t_H})} \ast
\rho_{t_F,t_H})(k) =
n_{\varepsilon,t_F,t_H}((S_{\gamma^*})_{t_F,\varepsilon
t_H},\gamma^*).
\end{multline*}
The second inequality follows of the proven one applied to the
domain $B\setminus S$, where $B$ is a sufficiently large ball,
containing $S$.
\end{proof}

For any $f\in {\mathcal S}(V^\bot)$, define its Fourier transform
$\hat{f}\in {\mathcal S}(V^\bot)$ by
\[
\hat{f}(\xi)=\int_{V^\bot}e^{-2\pi i(\xi,x)}f(x)\,dx.
\]

Recall the Poisson summation formula
\begin{equation}\label{e:Poisson}
\sum_{k\in \Gamma^\bot} f(k)=\frac{1}{|Q|}\sum_{k^*\in
{\Gamma^\bot}^*}\hat{f}(k^*), \quad f\in {\mathcal S}(V^\bot),
\end{equation}
where ${\Gamma^\bot}^*\subset V^{\bot}$ is the dual lattice of
$\Gamma^\bot $, and we used the relation \eqref{e:vol}.

Let us apply \eqref{e:Poisson} to the function
\begin{equation}\label{e:f}
f(x)=(\chi_{T_\varepsilon((S_{\gamma^*})_{t_F,\varepsilon t_H})}
\ast \rho_{t_F,t_H})(k_{\gamma^*}+x), \quad x\in V^\bot.
\end{equation}
The formula \eqref{e:Poisson} can be applied, because, for any
$N>0$, we have the estimate
\begin{equation}\label{e:Frho}
|\hat{\rho}_{t_F,t_H}(\xi)|\leq
C_N\frac{1}{1+t_F^N|\xi_F|^N+t_H^N|\xi_H|^N}, \quad \xi \in V^\bot.
\end{equation}

One has the relation
\[
\hat\chi_{T_\varepsilon((S_{\gamma^*})_{t_F,\varepsilon t_H})}(\xi)=
\varepsilon^{-q}e^{2\pi i(\xi,(1-T_{\varepsilon})(k_{\gamma^*}))}
\hat\chi_{(S_{\gamma^*})_{t_F,\varepsilon
t_H}}(T_{\varepsilon}(\xi)).
\]
Indeed, for any subset $S\subset V^\bot$, we have
\begin{align*}
\hat\chi_{T_\varepsilon(S)}(\xi) & =\int_{V^\bot}e^{-2\pi
i(\xi,x)}\chi_{T_\varepsilon(S)}(k_{\gamma^*}+x)dx\\
& =\int_{V^\bot}e^{-2\pi
i(\xi,x)} \chi_{S}(k_{\gamma^*}+(T_{\varepsilon^{-1}}(k_{\gamma^*})-k_{\gamma^*})+T_{\varepsilon^{-1}}(x))dx\\
& =\varepsilon^{-q}\int_{V^\bot}e^{-2\pi
i(\xi,T_\varepsilon(x^\prime)+(T_\varepsilon(k_{\gamma^*})-k_{\gamma^*}))}
\chi_{S}(k_{\gamma^*}+x^\prime)dx\\ &= \varepsilon^{-q} e^{2\pi
i(\xi,(1-T_{\varepsilon})(k_{\gamma^*}))}
\hat\chi_{S}(T_{\varepsilon}(\xi)).
\end{align*}
We also have
\[
\hat\rho_{t_F,t_H}(\xi)=\hat\rho_{1,1}(t_F\xi_F+t_H\xi_H), \quad
\xi\in V^\bot.
\]
By these relations, the Poisson formula \eqref{e:Poisson} applied to
the function $f$ given by \eqref{e:f} is written as
\begin{multline}\label{e:ne1}
n_{\varepsilon,t_F,t_H}((S_{\gamma^*})_{t_F,\varepsilon
t_H},\gamma^*) \\ = \frac{\varepsilon^{-q}}{|Q|} \sum_{k\in
{\Gamma^\bot}^*} e^{2\pi i(k^*,(1-T_{\varepsilon})(k_{\gamma^*}))}
\hat\chi_{(S_{\gamma^*})_{\varepsilon,t_F,t_H}}(
T_{\varepsilon}(k))\hat\rho_{1,1}(t_Fk_F+t_Hk_H).
\end{multline}
The series in the right hand side of~\eqref{e:ne1} converges
uniformly by the estimate~\eqref{e:Frho}.

Let us write
\[
n_{\varepsilon,t_F,t_H}((S_{\gamma^*})_{t_F,\varepsilon
t_H},\gamma^*)=n^{\prime}_{\varepsilon,t_F,t_H}
((S_{\gamma^*})_{t_F,\varepsilon
t_H},\gamma^*)+n^{\prime\prime}_{\varepsilon,t_F,t_H}
((S_{\gamma^*})_{t_F,\varepsilon t_H},\gamma^*),
\]
where
\begin{multline*}
n^{\prime}_{\varepsilon,t_F,t_H} ((S_{\gamma^*})_{t_F,\varepsilon
t_H},\gamma^*)\\ = \frac{\varepsilon^{-q}}{|Q|} \sum_{k\in
{\Gamma^\bot}^*, k_H= 0}\hat\chi_{(S_{\gamma^*})_{t_F,\varepsilon
t_H}}(T_{\varepsilon}(k))\hat\rho_{1,1}(t_Fk_F+t_Hk_H),
\end{multline*}
and
\begin{multline*}
n^{\prime\prime}_{\varepsilon,t_F,t_H}
((S_{\gamma^*})_{t_F,\varepsilon t_H},\gamma^*)\\ =
\frac{\varepsilon^{-q}}{|Q|} \sum_{k\in {\Gamma^\bot}^*, k_H\neq 0}
e^{2\pi i(1-\varepsilon^{-1})
(k_H,k_{\gamma^*})}\hat\chi_{(S_{\gamma^*})_{t_F,\varepsilon t_H}}(
T_{\varepsilon}(k)) \hat\rho_{1,1}(t_Fk_F+t_Hk_H).
\end{multline*}

Let $k\in {\Gamma^\bot}^*$ be such that $k_H=0$. Then $k\in F_V$.
Since ${\Gamma^\bot}^*\subset {\mathbb Q}^n$ and $F_V \cap {\mathbb
Q}^n=\{0\}$, we get $k=0$. Thus, we have
\begin{multline*}
n^{\prime}_{\varepsilon,t_F,t_H} ((S_{\gamma^*})_{t_F,\varepsilon
t_H},\gamma^*)\\
\begin{aligned}
 & =
\frac{\varepsilon^{-q}}{|Q|} \hat\chi_{(S_{\gamma^*})_{t_F,\varepsilon t_H}} (0)\\
& = \frac{\varepsilon^{-q}}{|Q|}{\rm vol}_{n-r}
((S_{\gamma^*})_{t_F,\varepsilon t_H})\\
& = \frac{\varepsilon^{-q}}{|Q|}{\rm vol}_{n-r} (P_{\gamma^*}\cap
S)+\frac{\varepsilon^{-q}}{|Q|}{\rm vol}_{n-r}
((S_{\gamma^*})_{t_F,\varepsilon t_H} \setminus S_{\gamma^*}).
\end{aligned}
\end{multline*}
We have the estimate 
\[
{\rm vol}_{n-r} ((S_{\gamma^*})_{t_F,\varepsilon t_H} \setminus
S_{\gamma^*})\leq C (t_F+t_H\varepsilon),
\]
therefore, we obtain that
\begin{equation}\label{e:n11}
n^{\prime}_{\varepsilon,t_F,t_H} ((S_{\gamma^*})_{t_F,\varepsilon
t_H},\gamma^*) = \frac{\varepsilon^{-q}}{|Q|}{\rm vol}_{n-r}
(P_{\gamma^*}\cap S)+O(t_F\varepsilon^{-q}+t_H\varepsilon^{1-q}).
\end{equation}

Consider the case when $k\in {\Gamma^\bot}^*$ and $k_H\neq 0$. For
any $t\in F_V$ and for any domain $D\subset P_{\gamma^*}$, we denote
\[
D(t)=\{x_H\in H : k_{\gamma^*}+t+x_H\in D\}\subset H.
\]
For any function $\phi\in {\mathcal S}(H)$, denote by $F_H(\phi)\in
{\mathcal S}(H)$ its Fourier transform:
\[
[F_H(\phi)](\xi_H)=\int_H \phi(x_H) e^{-2\pi i(\xi_H,x_H)}\,dx_H,
\quad \xi_H\in H.
\]
It is easy to see that
\[
[F_H(\chi_{D(t)})](\xi_H)=\int_H \chi_D(k_{\gamma^*}+t+x_H) e^{-2\pi
i(\xi_H,x_H)}\,dx_H, \quad \xi_H\in H.
\]
Therefore, we get
\begin{multline*}
\hat\chi_{(S_{\gamma^*})_{t_F,\varepsilon
t_H}}(T_{\varepsilon}(k))\\
\begin{aligned}
& =\int_{V^\bot}\chi_{(S_{\gamma^*})_{t_F,\varepsilon
t_H}}(k_{\gamma^*}+x)e^{-2\pi i(T_{\varepsilon}(k),x)}dx \\
& =\int_{F_V}\int_H \chi_{(S_{\gamma^*})_{t_F,\varepsilon
t_H}}(k_{\gamma^*}+x_F+x_H)e^{-2\pi i((k_F,x_F) +\varepsilon^{-1}(k_H,x_H))}dx_Fdx_H \\
& =\int_{F_V} e^{-2\pi i(k_F,x_F)}
F_H[\chi_{(S_{\gamma^*})_{t_F,\varepsilon t_H}( x_F)}]
(\varepsilon^{-1}k_H) dx_F,
\end{aligned}
\end{multline*}
and, as a consequence, we obtain the estimate
\begin{equation}\label{e:chi}
|\hat\chi_{(S_{\gamma^*})_{t_F,\varepsilon
t_H}}(T_{\varepsilon}(k))|\leq \int_{F_V}
|F_H[\chi_{(S_{\gamma^*})_{t_F,\varepsilon t_H}( x_F)}]
(\varepsilon^{-1}k_H)| dx_F.
\end{equation}
Hence, our considerations are reduced to the sufficiently well
studied problem of estimating the Fourier transform of the indicator
of a domain, and we can apply existing results.

\begin{proof}[Proof of Theorem~\ref{mainthm1}]
For any sufficiently small $\varepsilon>0$, $t_F>0$ and $t_H>0$, the
domain $(S_{\gamma^*})_{t_F,\varepsilon t_H}(x_F)$ has smooth
boundary. The Stokes formula allows us to write the Fourier
transform $F_H[\chi_{(S_{\gamma^*})_{t_F,\varepsilon t_H}(x_F)}]$ as
an oscillating integral over the boundary of
$(S_{\gamma^*})_{t_F,\varepsilon t_H}(x_F)$:
\begin{multline}\label{e:osc}
F_H[\chi_{(S_{\gamma^*})_{t_F,\varepsilon t_H}(x_F)}](\tau \omega)\\
=\frac{1}{\tau}\int_{\partial [(S_{\gamma^*})_{t_F,\varepsilon
t_H}(x_F)]} e^{-i\tau (\omega,x)} i_\omega (dx_1\wedge \ldots \wedge
dx_q),
\end{multline}
that implies the estimate
\begin{equation}\label{e:FH}
|F_H[\chi_{(S_{\gamma^*})_{t_F,\varepsilon t_H}(x_F)}](\xi)| =
O(|\xi|^{-1}), \quad |\xi|\to \infty.
\end{equation}
Hence, using the estimates \eqref{e:chi}, \eqref{e:FH} and
\eqref{e:Frho}, we get
\begin{multline}\label{e:n21}
|n^{\prime\prime}_{\varepsilon,t_F,t_H}
((S_{\gamma^*})_{t_F,\varepsilon
t_H},\gamma^*)|\\
\begin{aligned}
 & \leq  C \varepsilon^{-q} \sum_{k\in
{\Gamma^\bot}^*, k_H\neq 0}\varepsilon
|k_H|^{-1}\frac{1}{1+t_F^N|k_F|^N+t_H^N|k_H|^N}\\
& \leq  C \varepsilon^{1-q}
\int_{V^{\bot}} |x_H|^{-1} \frac{dx_F\,dx_H}{1+t_F^N|x_F|^N+t_H^N|x_H|^N}\\
& \leq C \varepsilon^{1-q} t_F^{-(p-r)}t_H^{1-q}.
\end{aligned}
\end{multline}
Putting $t_H=c_H$ to be a constant, independent of $\varepsilon>0$,
and $t_F=c_F\varepsilon^\alpha$, where
\[
\alpha=\frac{1}{p-r+1},
\]
by \eqref{e:n11} and \eqref{e:n21}, we obtain the statement of the
theorem.
\end{proof}

\begin{proof}[Proof of Theorem~\ref{mainthm}]
By assumption, the domain $S_{\gamma^*}(x_F)=S\cap
\{\gamma^*+x_F+H\}$ is strictly convex. Therefore, for any
sufficiently small $\varepsilon>0$, $t_F>0$ and $t_H>0$, the domain
$(S_{\gamma^*})_{t_F,\varepsilon t_H}(x_F)$ is strictly convex. By
the stationary phase method, we derive from the
representation~\eqref{e:osc} the following estimate
\begin{equation}\label{e:FH2}
|F_H[\chi_{(S_{\gamma^*})_{t_F,\varepsilon t_H}(x_F)}](\xi)| =
O(|\xi|^{-(q+1)/2}), \quad |\xi|\to \infty.
\end{equation}
Thus, using the estimate \eqref{e:chi}, \eqref{e:FH2} and
\eqref{e:Frho}, we obtain that
\begin{multline}\label{e:n22}
|n^{\prime\prime}_{\varepsilon,t_F,t_H} ((S_{\gamma^*})_{t_F,\varepsilon t_H},\gamma^*)|\\
\begin{aligned} & \leq  C \varepsilon^{-q} \sum_{k\in {\Gamma^\bot}^*, k_H\neq
0}\varepsilon^{(q+1)/2} |k_H|^{-(q+1)/2} \frac{1}{1+t_F^N|k_F|^N+t_H^N|k_H|^N}\\
& \leq  C \varepsilon^{-q} \varepsilon^{(q+1)/2}\int_{V^{\bot}}
|x_H|^{-(q+1)/2}\frac{dx_F\,dx_H}{1+t_F^N|x_F|^N+t_H^N|x_H|^N}\\
& \leq C \varepsilon^{-(q-1)/2} t_F^{-(p-r)}t_H^{-(q-1)/2}.
\end{aligned}
\end{multline}
Put $t_F=\varepsilon^{\alpha_F}$, $t_H=\varepsilon^{\alpha_H}$,
where $\alpha_F\geq 0$ and $\alpha_H\geq 0$ are chosen as follows.
If $\frac{q-1}{2}>p-r$, then
\[
\alpha_F=\frac{2q}{q+1+2(p-r)}, \quad
\alpha_H=\frac{q-1-2(p-r)}{q+1+2(p-r)}.
\]
If $\frac{q-1}{2}\leq p-r$, then
\[
\alpha_F=\frac{q+1}{2(p-r+1)}, \quad \alpha_H= 0.
\]
Using the estimates \eqref{e:n11} and \eqref{e:n22}, we immediately
conclude the proof.
\end{proof}

\bibliographystyle{amsalpha}

%%%%%%%%%%%%%%%%%%%%%%%%%%%%%%%%%%%%%%%%%%%%%%%%%%%%%%%%%%%%

\end{document}